\documentclass{conm-p-l}
\usepackage{amssymb}
\usepackage{amsxtra}
\usepackage[mathscr]{eucal}
\usepackage{enumerate}

\theoremstyle{plain}
\newtheorem{theorem}{Theorem}[section]
\newtheorem{lemma}[theorem]{Lemma}

\newtheorem{conjecture}[theorem]{Conjecture}

\theoremstyle{remark}
\newtheorem{remark}[theorem]{Remark}
\newtheorem{definition}[theorem]{Definition}
\newtheorem{example}[theorem]{Example}
\newtheorem*{remark*}{Remark}

\newtheorem{observation}[theorem]{Observation}

\numberwithin{equation}{section}

\input xy \xyoption{matrix} \xyoption{arrow}
          \xyoption{curve}  \xyoption{frame}

\newcommand\edge{\ar@{-}}
\newcommand\dashedge{\ar@{--}}

\newcommand\PP{{\mathbb P}}

\newcommand\A{{\mathcal A}}

\newcommand\calO{{\mathcal O}}

\newcommand\bfq{\mathbf{q}}

\newcommand\Oq{\calO_q}

\newcommand\kx{k^\times}

\newcommand\prim{\operatorname{Prim}}
\newcommand\spec{\operatorname{Spec}}

\newcommand\Pprim{\operatorname{Pprim}}
\newcommand\Pspec{\operatorname{Pspec}}
\newcommand\specP{\operatorname{Spec}_{\text{p}}}

\newcommand\chr{\operatorname{char}}

\newcommand\Aut{\operatorname{Aut}}

\newcommand\Obfqkn{\calO_{\bfq}(k^n)}
\newcommand\Obfqkxn{\calO_{\bfq}((\kx)^n)}

\newcommand\Oqktwo{\Oq(k^2)}
\newcommand\nrootone{\ne \root\bullet\of1}
\newcommand\SLtwok{SL_2(k)}
\newcommand\OqSLtwok{\Oq(SL_2(k))}
\newcommand\OqGLtwok{\Oq(GL_2(k))}

\newcommand\OqMtwok{\Oq(M_2(k))}

\newcommand\phihat{\widehat{\phi}}
\newcommand\Hspec{H\text{-}\spec}

\newcommand\Cl{{\rm{CL}}}
\newcommand\ftil{\widetilde{f}}
\newcommand\gtil{\widetilde{g}}
\newcommand\Vtil{\widetilde{V}}
\newcommand\id{\operatorname{id}}

\begin{document}

\title{Spectra of quantum algebras}

\author[K. R. Goodearl]{K. R. Goodearl}
\address{
Department of Mathematics \\
University of California\\
Santa Barbara, CA 93106 \\
U.S.A.
}
\email{goodearl@math.ucsb.edu}
\thanks{This research was supported
by US National Science Foundation grant DMS-1601184.}



\begin{abstract}
This is a survey of what is known and/or conjectured about the prime and primitive spectra of quantum algebras, of quantized coordinate rings in particular. The topological structure of these spectra, their relations to classical affine algebraic varieties, and their relations to each other are discussed.
\end{abstract}

\maketitle

\section{Introduction}

Suppose that an algebra $A$ is a quantized coordinate ring of an affine algebraic variety $V$ over an algebraically closed field. The prime and primitive spectra of $A$ are equipped with natural Zariski topologies, and the question arises of how these spaces are related to the corresponding classical spaces, namely the prime spectrum of the classical coordinate ring $\calO(V)$ (i.e., the space underlying the scheme structure of $V$) and the maximal ideal spectrum of $\calO(V)$ (which is homeomorphic to $V$ itself). We survey results and conjectures in several directions:
\begin{enumerate}
\item Are the spectra of $A$ topological quotients of the spectra of $\calO(V)$?
\item Are the spectra of $A$ homeomorphic to Poisson spectra of $\calO(V)$, where $\calO(V)$ is equipped with a Poisson structure from a semiclassical limit process?
\item Can the topologies on the spectra of $A$ be described via piecewise classical data, from finite collections of classical spectra and connections among them?
\end{enumerate}

Positive answers to (1) are conjectured in general, while for (2), positive answers are only conjectured for generic cases in characteristic zero. Under an additional hypothesis that holds for many quantized coordinate rings, any positive answer to (2) implies a positive answer to (1). Positive answers to (3), finally, are conjectured for generic cases in arbitrary characteristic. In characteristic zero, a positive answer to (3), coupled with a matching positive answer to the corresponding question about the Poisson spectra of $\calO(V)$, implies a positive answer to (2).

\subsection{}
Thoughout the paper, we work over an algebraically closed base field $k$.

\subsection{}
If $R$ is the coordinate ring of a classical affine variety $V$ over $k$, then the space of maximal ideals of $R$ is homeomorphic to $V$, and the space of prime ideals of $R$ is homeomorphic to the scheme of irreducible subvarieties of $V$. In the noncommutative world, primitive ideals, rather than just maximal ideals, are the appropriate ideals to relate to points, so it is the primitive spectrum rather than the maximal ideal spectrum which takes over the role of the variety.

For a noncommutative algebra $A$, set
\begin{align*}
\prim A &:= \{\; \text{primitive ideals of}\; A\;\}, 
&\spec A &:= \{\; \text{prime ideals of}\; A\;\},
\end{align*}
both equipped with the natural Zariski topologies. (In case there is asymmetry of primitivity, let us interpret ``primitive" as ``left primitive".)

\section{Three examples}  \label{expls}

To illustrate prototypical patterns, we exhibit the spectra of three basic examples -- the quantum plane, quantum $SL_2$, and quantum $GL_2$, meaning quantized coordinate rings of the plane $k^2$ and the algebraic groups $SL_2(k)$ and $GL_2(k)$. We just give the standard (single parameter) versions, with parameter $q \in \kx$.

\begin{example} \label{exqplane}
The quantum plane is the algebra
$$A = \Oqktwo := k\langle x,y \mid xy = qyx \rangle.$$
In the \emph{generic} case, i.e., if $q$ is not a root of unity, $\spec A$ may be drawn as in Figure \ref{specOqk2}, where dashed lines indicate inclusions (e.g., \cite[Example II.1.2]{KBKG02}).
\begin{figure}[ht]
$$\xymatrixrowsep{4pc}\xymatrixcolsep{2.5pc}
\xymatrix{
 \edge[r] &\langle x,\, y-\beta \rangle \save+<0ex,-4ex>
\drop{(\beta \in \kx)} \restore \edge[r] &&\langle x,y \rangle
&\edge[r] &\langle x-\alpha,\, y \rangle
\save+<0ex,-4ex> \drop{(\alpha \in \kx)} \restore
\edge[r] & \\
&\langle x\rangle \dashedge[ul] \dashedge[ur] \dashedge[urr] &&&&\langle y\rangle
\dashedge[ull] \dashedge[ul] \dashedge[ur]\\
 &&&\langle 0\rangle \dashedge[ull] \dashedge[urr]
}$$
\caption{$\spec \Oqktwo$}
\label{specOqk2}
\end{figure}

The primitive spectrum of $A$ consists of the maximal ideals together with the zero ideal (e.g., \cite[Example II.7.2]{KBKG02}). We picture this spectrum as in Figure \ref{primOqk2} -- the $x$- and $y$-axes of the classical plane, together with a fat point, which is dense.
\begin{figure}[ht]
$$\xymatrixrowsep{0.5pc}\xymatrixcolsep{2.5pc}
\xymatrix{
\edge[ddddrrrrrr] &&&&&&\edge[ddddllllll]  \\  \\  \\  \\  
&&&&&& \\  \\
&&&\bullet 
}$$
\caption{$\prim \Oqktwo$}
\label{primOqk2}
\end{figure}
There is an obvious surjection $\pi : k^2 \rightarrow \prim A$, where $\pi(\alpha,0) = \langle x-\alpha,\, y \rangle$ and $\pi(0,\beta) = \langle x,\, y-\beta \rangle$ for $\alpha,\beta \in k$ while $\pi(\alpha,\beta) = \langle 0 \rangle$ for $(\alpha,\beta) \in (\kx)^2$. It is easily checked that $\pi$ is a topological quotient map. In other words, the topology on $\prim \Oqktwo$ is the quotient topology that results from $\pi$.
\end{example}

Quantum $SL_2$ and quantum $GL_2$ are obtained from the quantum $2\times2$ matrix algebra $\OqMtwok$, which we present via generators $a$, $b$, $c$, $d$ and the following relations:
\begin{align*}
ab &= qba  &ac &= qca  &bc &= cb  \\
bd &= qdb  &cd &= qdc  &ad-da &= (q-q^{-1})bc.
\end{align*}

The element $D_q := ad-qbc$ in this algebra is a central element known as the ($2\times2$) \emph{quantum determinant}. In parallel with the classical commutative case, quantized coordinate rings of $SL_2(k)$ and $GL_2(k)$ are defined as follows:
\begin{align*}
\OqSLtwok &:= \OqMtwok/\langle D_q-1 \rangle  &\OqGLtwok &:= \OqMtwok[D_q^{-1}].
\end{align*}
It is convenient to use the same symbols $a$, $b$, $c$, $d$ for the cosets $a+\langle D_q-1 \rangle,\dots$ in $\OqSLtwok$.

\begin{example} \label{exoqsl2}
Consider $A = \OqSLtwok$ in the generic case, $q\nrootone$. We may picture $\spec A$ as in Figure \ref{specOqSL2}, where the dashed lines again indicate inclusions (e.g., \cite[Example II.1.3]{KBKG02}).
\begin{figure}[ht]
$$\xymatrixrowsep{2.5pc}\xymatrixcolsep{2pc}
\xymatrix{
\edge[rr] &&\langle a-\lambda,\,b,\, c,\,
d-\lambda^{-1} \rangle \save+<0ex,-4ex> \drop{(\lambda \in
k^\times)} \restore \edge[rr] && \\
 &&\langle b,c \rangle \dashedge[ull] \dashedge[urr]\\
\dashedge[urr] \edge[rr] &&\langle \gamma b-
\beta c \rangle \save+<0ex,-4ex> \drop{([\beta:\gamma] \in \PP^1(k))}
\restore \edge[rr] &&\dashedge[ull]\\
 &&\langle 0\rangle \dashedge[ull] \dashedge[urr]
}$$
\caption{$\spec \OqSLtwok$}
\label{specOqSL2}
\end{figure}
Here the primitive ideals consist of the maximal ideals together with the height $1$ prime ideals $\langle \gamma b-
\beta c \rangle$ (e.g., \cite[Example II.8.6]{KBKG02}). $\prim A$ is pictured in Figure \ref{primOqSL2}.
\begin{figure}[ht]
$$\xymatrixcolsep{2pc} \xymatrixrowsep{8pc}
\xymatrix{
\edge[rrrr] &&&&\save+<0ex,0ex> \drop{\circ} \restore \edge[rrrr] &&&&  \\
\infty \edge[rrrrrrrr] &&\save+<0ex,0ex> \drop{\bullet} \restore \dashedge[ull]\dashedge[urrrrrr] &&\save+<0ex,0ex> \drop{\bullet} \restore \dashedge[ullll]\dashedge[urrrr] &&\save+<0ex,0ex> \drop{\bullet} \dashedge[ullllll]\dashedge[urr] \restore &&\infty
}$$
\caption{$\prim \OqSLtwok$}
\label{primOqSL2}
\end{figure}
The punctured line at the top of the figure represents $k \setminus \{0\}$, while the bottom line stands for $\PP^1(k)$.  The bullets are sample points, and the inclusions represented by the dashed lines indicate that the closure of each point on the projective line consists of that point together with the full punctured affine line.

As in Example \ref{exqplane}, there is a topological quotient map $\pi : SL_2(k) \rightarrow \prim A$, this time given by 
\begin{equation}  \label{tquoSLtwok}
\pi \left(\left[\begin{smallmatrix} \alpha&\beta\\ \gamma&\delta \end{smallmatrix}\right]\right) = \begin{cases}
\langle a-\alpha,\,b,\,c,\, d-\delta \rangle &\quad(\beta= \gamma= 0)\\
\langle \gamma b- \beta c \rangle &\quad(\beta,\gamma\; \text{not both}\;=0).
\end{cases}
\end{equation}
\end{example}

\begin{example}  \label{exoqgl2}
Finally, consider $A = \OqGLtwok$, $q\nrootone$. We display $\spec A$ without indicating all inclusions. As a set, it is the disjoint union of $4$ subsets as drawn in Figure \ref{specOqGL2} (e.g., \cite[Example II.8.7]{KBKG02}).
\begin{figure}[ht]
$$\xymatrixrowsep{2.5pc}\xymatrixcolsep{1.8pc}
\xymatrix{
\edge[rr] &&\langle a-\alpha,\,b,\, c,\,
d-\delta \rangle \save+<6.5ex,-4ex> \drop{(\alpha,\delta \in
k^\times)} \restore \edge[rr] && &&\edge[r] &\langle b, ad-\lambda \rangle \save+<0ex,-4ex> \drop{(\lambda \in k^\times)} \restore \edge[r] & \\
&\edge[r] &\langle b,c,p_{\alpha,\delta}\rangle \save+<0ex,-4ex> \drop{(\alpha,\delta \in
k^\times)} \restore \dashedge[u] \edge[r] & &&& &\langle b \rangle \dashedge[ul] \dashedge[ur]  \\
&&\langle b,c \rangle \dashedge[ul] \dashedge[ur]  \\
\edge[rr] &&\langle b-\mu c,\,D_q-\lambda \rangle \save+<6.5ex,-4ex> \drop{(\mu,\lambda \in
k^\times)} \restore \edge[rr] && &&\edge[r] &\langle c, ad-\lambda \rangle \save+<0ex,-4ex> \drop{(\lambda \in k^\times)} \restore \edge[r] & \\
&\edge[r] &\langle p'_{\mu,\lambda}\rangle \save+<0ex,-4ex> \drop{(\mu,\lambda \in
k^\times)} \restore \dashedge[u] \edge[r] & &&& &\langle c \rangle \dashedge[ul] \dashedge[ur]  \\
&&\langle 0 \rangle \dashedge[ul] \dashedge[ur]
}$$
\caption{$\spec \OqGLtwok$}
\label{specOqGL2}
\end{figure}
In the upper left subdiagram, $p_{\alpha,\delta}$ denotes an irreducible polynomial in $k[a,d]$
such that $p_{\alpha,\delta}(\alpha,\delta) = 0$,
while in the lower left subdiagram, $p'_{\mu,\lambda}$ denotes a polynomial in $k[b,c,D_q]$ of the form $f(bc^{-1},D_q)c^m$ where $f$ is an irreducible polynomial in $k[s,t]$ of $s$-degree $m$ with $f(\mu,\lambda) = 0$.

In this algebra, the primitive ideals are the prime ideals appearing at the top levels of the $4$ diagrams in Figure \ref{specOqGL2} (e.g., \cite[Example II.8.7]{KBKG02}; cf.~Theorem \ref{stratif}). The primes $\langle a-\alpha,\,b,\, c,\,
d-\delta \rangle$ are maximal ideals, while the remaining primitive ideals are primes of height $2$. The latter may be grouped together as
$$\bigl\{ \langle \gamma b - \beta c,\, D_q-\lambda \rangle \bigm| [\beta:\gamma] \in \PP^1(k),\, \lambda \in \kx \bigr\}.$$
The only inclusions between primitive ideals in different subsets of $\prim A$ are 
$$\langle \gamma b - \beta c,\, D_q-\lambda \rangle \subseteq \langle a-\alpha,\,b,\, c,\,
d-\alpha^{-1}\lambda \rangle \qquad \forall\ [\beta:\gamma] \in \PP^1(k),\ \alpha,\lambda \in \kx.$$
This leads to a rough picture of $\prim A$ as in Figure \ref{primOqGL2}, where inclusions between primes in the lower and upper layers are not indicated.
\begin{figure}[ht]
$$\xymatrixrowsep{1pc}\xymatrixcolsep{1pc}
\xymatrix{
&&&&&& \edge[dlllll] \edge[dddl]  \\
&\edge[dddl]  \\
&&&(\kx)^2  \\
&&&\ar@{{o}{o}{o}}[ddd] &&\edge[dlllll]  \\
&  \\
&&&&&& \edge[dlllll] \edge[dddl]  \\
&\edge[dddl] &&  \\
&&&\PP^1(k) \times \kx  \\
&&&&&\edge[dlllll]  \\
& 
}$$
\caption{$\prim \OqGLtwok$}
\label{primOqGL2}
\end{figure}

The topology of the space in this picture may be described in terms of the Zariski topologies on $\PP^1(k) \times \kx$ and $(\kx)^2$, together with the map
\begin{align*}
\phi : \{\,\text{closed subsets of}\; \PP^1(k) \times \kx\,\} &\longrightarrow \{\,\text{closed subsets of}\; (\kx)^2\,\}  \\
X &\longmapsto \bigcup_{( [\nu:\mu],\lambda ) \in X} \{ (x,y) \in (\kx)^2 \mid xy = \lambda \}.
\end{align*}
(That the union in the display is closed in $(\kx)^2$ follows from the fact that the projection $\PP^1(k) \times \kx \rightarrow \kx$ is a closed map.)
Namely, a set $X \sqcup Y$, where $X \subseteq \PP^1(k) \times \kx$ and $Y \subseteq (\kx)^2$, is closed in $\prim A$ if and only if
\begin{itemize}
\item $X$ is closed in $\PP^1(k) \times \kx$;
\item $Y$ is closed in $(\kx)^2$;
\item $\phi(X) \subseteq Y$.
\end{itemize}
See Section \ref{piecewise} for a description of a general framework into which this example fits. 

There is a topological quotient map $GL_2(k) \rightarrow \prim A$, parallel to \eqref{tquoSLtwok}.
\end{example}

\section{Topological quotient properties}

As observed in Example \ref{exqplane}, the primitive ideal space of the quantum plane, $\prim \Oqktwo$, is a topological quotient of the classical plane, $k^2$, for $q\nrootone$.
Modulo a minor technical condition, this quotient picture carries over to quantum affine spaces in general, including the multiparameter algebras 
\begin{equation}  \label{Obfqkn}
\Obfqkn := k \langle x_1,\dots,x_n \mid x_i x_j = q_{ij} x_j x_i\ \forall\ i,j \rangle,
\end{equation} 
where $\bfq = (q_{ij})$ is a multiplicatively skew-symmetric matrix in $M_n(\kx)$.

\begin{theorem} \label{topoquo1} \textup{\cite[Theorems 4.11, 6.3]{KGEL00b}}
Let $A = \Obfqkn$, and assume either $-1 \notin \langle q_{ij} \rangle$ or $\chr k = 2$.  Then $\prim A$ is a topological quotient of $k^n$ and $\spec A$ is a {\rm(}compatible\/{\rm)} topological quotient of $\spec \calO(k^n)$. Compatibility means that there is a topological quotient map $\pi : \spec \calO(k^n) \rightarrow \spec A$ such that the restriction of $\pi$ to $\max \calO(k^n) \approx k^n$ is a topological quotient map onto $\prim A$.

These statements hold more generally for cocycle twists of commutative affine algebras graded by torsionfree abelian groups, and hence for quantum toric varieties.
\end{theorem}

The technical condition in Theorem \ref{topoquo1} excludes elements of order $2$ from the subgroup $\langle q_{ij} \rangle$ of $\kx$. It is not known whether this condition is necessary. 

Some other cases following this pattern are known -- e.g., the primitive spectrum of $\OqSLtwok$, $q\nrootone$, is a topological quotient of the variety $\SLtwok$ \cite[Example 1.1]{KG01}.

We conjecture that such topological quotient results in fact hold much more widely among quantized coordinate rings:

\begin{conjecture} \label{quoconj}  \textup{(\cite[Introduction]{KGEL00b}, \cite[Conjecture 1.6]{KG01})}
Let $A$ be a quantized coordinate ring of an affine variety $V$. Then $\spec A$ and $\prim A$ are compatible topological quotients of $\spec \calO(V)$ and $\max\calO(V) \approx V$.
\end{conjecture}

At this point, there is no accepted or even proposed definition of what it should mean for an algebra to be counted as a quantized coordinate ring. To clarify: the conjecture above is  to apply at least to currently accepted families of quantized coordinate rings, such as quantum matrices, quantum semisimple groups and quantum Borel subgroups thereof, as well as quantum symplectic and euclidean spaces.

Beyond Theorem \ref{topoquo1}, Conjecture \ref{quoconj} has been confirmed for quantum symplectic and euclidean spaces in \cite[Corollaries 5.9, 5.10]{Oh08}, \cite[Corollary 10]{OhMPa}. 
A related conjecture which we discuss below (Conjecture \ref{semiclassconj}) leads to additional cases. More precisely, whenever Conjecture \ref{semiclassconj} and an additional hypothesis hold, Conjecture \ref{quoconj} follows (see Remark \ref{semiclass.implies.topoquo}). Confirmed cases include $\Oq(SL_2(k))$, $\Oq(GL_2(k))$ and $\Oq(SL_3(k))$, as noted below.

\section{Poisson relations}  \label{semiclass.Poisson}

One approach to confirming Conjecture \ref{quoconj} is to identify a space $X_s(V)$ and a subspace $X_p(V)$ such that (1) $X_s(V)$ and $X_p(V)$ are compatible topological quotients of $\spec \calO(V)$ and $\max\calO(V)$; and (2) there are compatible homeomorphisms $X_s(V) \approx \spec A$ and $X_p(V) \approx \prim A$. Knowledge of the structures of $X_s(V)$ and $X_p(V)$ then also provides knowledge about the structures of $\spec A$ and $\prim A$.

The main instance of this approach will be outlined in the present section. It involves (a) constructing a \emph{semiclassical limit} of a family of quantum algebras; (b) imposing a \emph{Poisson structure} that arises from the construction process; (c) investigating \emph{Poisson-prime} and \emph{Poisson-primitive spectra} of such a semiclassical limit; (d) relating the spaces obtained in (c) to prime and primitive spectra of \emph{generic} members of the original family. Although semiclassical limits and their Poisson structures can be obtained for general families, the relations implementing (d) do not work well in positive characteristic. Further, clarification of the attribute ``generic" in (d) is required. We address these concepts and points separately.

\subsection{Poisson algebras}  \label{Poisson.alg}
\begin{definition}  \label{poissonthings}
 A \emph{Poisson algebra} is an algebra $R$ equipped with a Lie algebra bracket $\{-,-\} : R \times R \rightarrow R$ such that
$$\{a,xy\} = \{a,x\}y + x\{a,y\}\ \ \forall\; a,x,y \in R.$$
Due to the antisymmetry of the bracket, $\{xy,a\} = x\{y,a\} + \{x,a\}y$ for all $a,x,y \in R$ as well.

A \emph{Poisson ideal} in $R$ is any ideal $I$ such that $\{R,I\} \subseteq I$.

A \emph{Poisson-prime ideal} in $R$ is any proper Poisson ideal $P$ such that 
$$IJ \subseteq P \implies I \subseteq P\ \text{or}\ J \subseteq P,\ \ \forall\ \text{Poisson ideals}\ I,J \subseteq R.$$
In case $R$ is noetherian and $\chr k = 0$, the Poisson-prime ideals of $R$ are precisely those Poisson ideals which are also prime ideals (e.g., \cite[Lemma 1.1]{KG06}).

A \emph{Poisson-primitive ideal} of $R$ is  any ideal $P$ which is the largest of the Poisson ideals contained in some maximal ideal of $R$. Such ideals are also Poisson-prime.

The \emph{Poisson-prime} and \emph{Poisson-primitive spectra} of $R$ are
\begin{align*}
\Pspec R &:= \{\; \text{Poisson-prime ideals of}\ R\;\}  \\
\Pprim R &:= \{\; \text{Poisson-primitive ideals of}\ R\;\} \subseteq \Pspec R.
\end{align*}
Both of these sets are equipped with Zariski topologies.
\end{definition}

These Poisson spectra form natural topological quotients of classical spectra in many cases:

\begin{theorem} \label{topoquo2} \textup{\cite[Theorems 4.1, 1.5]{KG06}, \cite[Lemma 9.3]{KG10}}
Let $R$ be a commutative noetherian Poisson $k$-algebra, with $\chr k = 0$. Then $\Pspec R$ is a topological quotient of $\spec R$, via the map
$$\pi : P \longmapsto (\; \text{\rm largest Poisson ideal}\; \subseteq P\;).$$

Now assume in addition that $R$ is affine over $k$ and satisfies the following portion of the Poisson Dixmier-Moeglin Equivalence: all $P \in \Pprim R$ are locally closed in $\Pspec R$. Then $\Pprim R$ is a topological quotient of $\max R$, via $\pi$.
\end{theorem}

See, e.g., \cite{BLLM, KBIG, KG06, KGSL11, DJSQO, SLOLS, LWW, Oh99, Oh06, Oh08}  for some situations where the Poisson Dixmier-Moeglin Equivalence has been established.

\subsection{Semiclassical limits}  \label{semiclass}
These ``limits" arise from parametrized families of algebras. We review the process of forming semiclassical limits in the one-parameter case (a discussion of multiparameter cases is given in \cite{KG10}).

\begin{definition} Let $\A$ be a torsionfree algebra over a Laurent polynomial ring $k[t^{\pm1}]$ such that $\A/(t-1)\A$ is commutative. The family $\bigl( \A_q := \A/(t-q)\A \bigr)_{q \in \kx}$ is the corresponding \emph{flat family} of $k$-algebras, and the commutative algebra $\A_1$ is the algebra underlying the \emph{semiclassical limit} of this family.

This semiclassical limit preserves a trace of the (possible) noncommutativity in the family $\bigl( \A_q \bigr)$, expressed via an induced Poisson bracket.
Namely, since all additive commutators $[a,b]$ ($= ab-ba$) in $\A$ are divisible by $t-1$, and uniquely so due to the torsionfreeness assumption on $\A$, there is a bilinear bracket $\frac1{t-1}[-,-]$ on $\A$. It is easily checked that $\frac1{t-1}[-,-]$ is a Poisson bracket on $\A$, and that it induces a Poisson bracket $\{-,-\}$ on $\A_1$.

More precisely, it is the Poisson algebra $\bigl( \A_1,\, \{-,-\} \bigr)$ which is the \emph{semiclassical limit} of the family $\bigl( \A_q \bigr)_{q \in \kx}$.
\end{definition}

\begin{example}
 The Poisson structure on $\calO(\SLtwok)$, which is the semiclassical limit of the family $\bigl( \OqSLtwok \bigr)_{q\in\kx}$ corresponding to the algebra $\A = \calO_t(SL_2(k[t^{\pm1}]))$, is determined by the following data:
\begin{align*}
\{a,b\} &= ab  &\qquad \{a,c\} &= ac &\qquad \{a,d\} &= 2 bc  \\
\{b,c\} &= 0 &\{b,d\} &= bd &\{c,d\} &= cd\,.  
\end{align*}
\end{example}

\subsection{The qualifier ``generic"}
\emph{Question}: Which members of a family of quantum algebras should be labelled ``generic"? To our knowledge, no rigorous answer to this question has been given. The underlying idea is to apply the label ``generic" to those members which exhibit the ``most general behavior" within the family -- in particular, the ``most general noncommutativity" (i.e., the ``least commutativity"). The choices that respect this principle are relatively obvious in many cases, much less so in others.

The most common usage of the term ``generic" occurs within the ``standard" families of quantized coordinate rings, families indexed by a single parameter. 
These include the examples displayed in Section \ref{expls}, which are members of the families $\bigl( \Oqktwo \bigr)_{q\in\kx}$, $\bigl( \OqSLtwok \bigr)_{q\in\kx}$ and $\bigl( \OqGLtwok \bigr)_{q\in\kx}$ respectively. 
Among the ``standard quantized coordinate rings" are the algebras named (for short) quantum affine spaces, quantum tori, quantum matrices, quantum semisimple groups, quantum general linear groups, quantum Borel subgroups of semisimple groups, quantum symplectic and euclidean spaces. (See, e.g., \cite{KBKG02, KG00} for surveys.) They are typically denoted $\Oq(V)$, $k_q[V]$ or $R_q[V]$, where $V$ is an affine variety such as $k^n$, $(\kx)^n$, $M_{m,n}(k)$, a semisimple algebraic group $G$ or a Borel subgroup $B^\pm$ of $G$, a general linear group $GL_n(k)$, etc.
Here is the common usage within these families:

\begin{definition}  \label{generic}
Let $\Oq(V)$ be a ``standard quantized coordinate ring". Within this setting, the parameter $q$ and the algebra $\Oq(V)$ are said to be \emph{generic} when $q$ is not a root of unity.
\end{definition}

For an immediate example, observe that in $\Oqktwo$, powers $x^i$ and $y^j$ of the given generators commute if and only if $q^{ij} = 1$. Thus, if $q$ is generic, $x^i$ and $y^j$ commute only when one of them equals $1$. In contrast, when $q$ is an $l$-th root of unity, $x^i$ and $y^j$ commute whenever $l$ divides $ij$.

The properties of the generic members of a family typically differ vastly from those of the non-generic members. For instance, the generic quantum planes $\Oqktwo$ have center $k$ (as is easily checked) and are primitive (e.g., \cite[Exercise 3ZA]{intronoeth}). Since these algebras are infinite dimensional over their centers, Kaplansky's Theorem (e.g., \cite[Theorem 13.3.8]{JMJR}) shows that they cannot satisfy any polynomial identity. On the other hand, all the non-generic quantum planes do satisfy polynomial identities, because these algebras are finitely generated modules over their centers. (If $q$ is a primitive $l$-th root of unity, the center of $\Oqktwo$ is $k[x^l,y^l]$.) 

For another concrete instance, consider the quantum tori of rank two,
$$\Oq((\kx)^2) := k \langle x^{\pm1}, y^{\pm1} \mid xy = qyx \rangle,$$
for $q\in \kx$. If $q$ is generic, $\Oq((\kx)^2)$ is a central simple $k$-algebra with Krull and global dimension $1$ (\cite[Corollary 1.5, Theorems 2.1, 2.4]{Jat}, \cite[Lemma 1.1]{Lor}, \cite[Proposition 1.3, Corollary 3.10]{JMJP}). If $q$ is not generic, one finds that $\Oq((\kx)^2)$ has infinitely many ideals, and its Krull and global dimensions are $2$ (e.g., \cite[Theorems 2.3, 2.4]{Jat}).

\begin{remark}
That the generic members of the standard families are the ones indexed by non-roots of unity is an artifact of the way the parameter $q$ is used in the definitions of these families. For example, in the standard quantum plane $\Oqktwo$, the scalar $q$ is the coefficient of the (non)commutativity constraint $xy=qyx$, but we can get a similar constraint on replacing $q$ by any scalar-valued function. For simplicity, suppose we build a family $(A_q)_{q\in k}$ where $A_q := k \langle x,y \mid xy = (3q-2) yx\rangle$. Here one should define the generic members of the family to be those $A_q$ such that $3q-2$ is nonzero and not a root of unity, rather than those for which $q$ is not a root of unity. Going to an extreme, one could even argue that for $\chr k \ne 3$, only the single member $A_{2/3}$ should be called generic, since it is a non-noetherian ring with zero-divisors, whereas every other $A_q$ is a noetherian domain.
\end{remark}

A number of multiparameter families of quantized coordinate rings occur in the literature (see, e.g., \cite{KBKG02, KG00, KGSL11, Hrt}). However,
 the question of which members of such families should be considered ``generic" has not even been addressed. For instance, consider  the multiparameter quantum affine spaces $\Obfqkn$ as in \eqref{Obfqkn} and the corresponding quantum tori
$$\Obfqkxn := k \langle x_1^{\pm1},\dots,x_n^{\pm1} \mid x_i x_j = q_{ij} x_j x_i\ \forall\ i,j \rangle.$$ 
At first glance, taking the one-parameter case as a guide, one might consider as ``generic" the algebras $\Obfqkn$ and $\Obfqkxn$ where none of the $q_{ij}$ for $i\ne j$ is a root of unity. However, this still leaves room for ``non-general" relations among the $q_{ij}$, such as $q_{12} = q_{13}^{-1}$, which lets $x_1$ commute with $x_2x_3$. If also $q_{12} = q_{13}^{-1} = q_{23}$, then $x_1x_2x_3$ is central in $\calO_\bfq(k^3)$ and $\calO_\bfq((\kx)^3)$.

The ``most general" multiplicatively skew-symmetric matrices $\bfq \in M_n(\kx)$ are those for which the abelian group $\langle q_{ij} \rangle \le \kx$ is free of rank $n(n-1)/2$ (in which case the $q_{ij}$ with $i<j$ form a basis). Let us tentatively call these $\bfq$, as well as the corresponding algebras $\Obfqkn$ and $\Obfqkxn$, \emph{generic}.

To buttress this choice, note that the non-generic quantum tori $\Obfqkxn$ always have ``more commutativity" than the generic ones. For example, let $n=3$ and assume $\bfq$ is not generic, whence $q_{12}^a q_{13}^b q_{23}^c = 1$ for some nonzero integer triple $(a,b,c)$. If $a=0$, then $x_1^b x_2^c$ commutes with $x_3$, and similar nontrivial relations occur if $b$ or $c$ is zero. In any case, $x_1^b x_2^c$ commutes with $x_2^a x_3^b$, which gives a nontrivial relation if at most one of $a$, $b$, $c$ is zero. In all these non-generic cases, the group of units of $\calO_\bfq((\kx)^3)$ contains a non-cyclic abelian subgroup $G$ such that $G \cap \kx = \{1\}$.

In the generic case, the center of $\Obfqkxn$ is $k$, and $\Obfqkxn$ is a simple algebra with Krull and global dimension $1$ \cite[Proposition 1.3, Corollary 3.10]{JMJP}, although these properties also extend to certain non-generic cases (e.g., \cite[\S3.11]{JMJP}).

\subsection{Primary semiclassical limit conjecture}  \label{semiclasslimconj}
A prerequisite for a useful confirmation or refutal of the following conjecture is to pin down a specification of the term ``generic". (For standard families of quantized coordinate rings, Definition \ref{generic} is to be used.) To emphasize this point, we adopt the phrase \emph{suitably generic} in formulating this and other conjectures.

\begin{conjecture} \label{semiclassconj}  \textup{\cite[Conjecture 9.1]{KG10}}
Assume that $\chr k = 0$. Let $A$ be a suitably generic member of a flat family of quantized coordinate rings for an affine variety $V$ over $k$, with semiclassical limit $\calO(V)$.  
Then there exist compatible homeomorphisms $\spec A \longrightarrow \Pspec \calO(V)$ and $\prim A \longrightarrow \Pprim \calO(V)$.
\end{conjecture}

We list a number of cases in which Conjecture \ref{semiclassconj} is known to hold:
\begin{itemize}
\item Quantum affine spaces and quantum affine toric varieties \cite[Theorem 3.5, Corollary 3.6]{OPS}, \cite[Theorems 3.2, 4.1, 5.1]{KGEL09};
\item Quantum symplectic and euclidean spaces \cite[Theorem 5.8]{Oh08}, \cite[Corollary 9]{OhMPa};
\item $\OqSLtwok$ and $\Oq(GL_2(K))$ \cite[Example 9.7, Evidence 9.8(a)]{KG10}, \cite[Corollary 1.3]{Fry}; 
\item $\Oq(SL_3(K))$ \cite[Theorem 1.2]{Fry}.
\end{itemize}

\begin{remark}  \label{semiclass.implies.topoquo}
Conjecture \ref{semiclassconj} implies Conjecture \ref{quoconj} in the following setting. Suppose $k$, $A$, $V$ are as in the above conjecture and that, additionally, $\calO(V)$ satisfies the Poisson Dixmier-Moeglin Equivalence. In view of Theorem \ref{topoquo2}, we find: If Conjecture \ref{semiclassconj} holds, then Conjecture \ref{quoconj} holds as well.
\end{remark}

\begin{remark}  \label{char2probs}
We briefly illustrate how Conjecture \ref{semiclassconj} may fail in positive characteristic.

(a) Let $A = \Oqktwo$, where $q \nrootone$ but $\chr k = 2$, and $R = \calO(k^2) = k[x,y]$. Due to the characteristic, the elements $x^2$ and $y^2$ in $R$ are Poisson-central, meaning that $\{x^2,R\} = \{y^2,R\} = 0$.

As one can check, the ideals $\langle x^2+\lambda^2y^2\rangle$ of $R$, for $\lambda \in \kx$, are Poisson-prime but not prime. Since $k$ is infinite, $\bigcap_{\lambda\in\kx} \langle x^2+\lambda^2y^2\rangle = \langle0\rangle$, from which it follows that $\{\langle0\rangle\}$ is not open in $\Pspec R$. 

In both $\spec A$ and $\Pspec R$, the zero ideal is the unique dense point. However, $\langle0\rangle$ is open in $\spec A$ (as is clear from Figure \ref{specOqk2}). Therefore $\spec A$ is not homeomorphic to $\Pspec R$ in this case.

(b) The example in (a) can be ``rescued" by replacing the space $\Pspec R$ by the space $\specP R$ of prime Poisson ideals of $R$. The set $\specP R$ consists of those prime ideals of $R$ which have generating sets matching those displayed in Figure \ref{specOqk2}, and so there is an obvious homeomorphism $\spec A \longrightarrow \specP R$.

(c) Keep $\chr k = 2$, but let $A$ be a generic member of the family $\bigl( \calO_{q^2}(k^2) \bigr)_{q\in\kx}$.  In this case, the semiclassical limit is $R = k[x,y]$ with zero Poisson bracket, so $\Pspec R = \spec R$, which is not homeomorphic to $\spec A$.
\end{remark}

\subsection{A secondary semiclassical limit conjecture} 
Whenever Conjecture \ref{semiclassconj} holds, all of the considered generic quantized coordinate rings of $V$ must have homeomorphic prime and primitive spectra. Thus, we are led to a related conjecture:

\begin{conjecture} \label{diffgenericconj}
Assume that $\chr k = 0$.  Let $A$ and $A'$ be suitably generic members of some flat family of quantized coordinate rings for an affine variety $V$ with semiclassical limit $\calO(V)$. Then there exist compatible homeomorphisms $\spec A \longrightarrow \spec A'$ and $\prim A \longrightarrow \prim A'$.
\end{conjecture}

Evidently, Conjecture \ref{diffgenericconj} holds in all the cases where Conjecture \ref{semiclassconj} has been verified.  Another positive scenario is known, although the proof involves a shift of the base field. In this case, we take ``suitably generic" to mean ``transcendental over a suitable subfield of $k$".

\begin{observation} \label{transmatch}
Let $\bigl( A_q \bigr)_{q\in\kx} = \bigl( \A/(t-q)\A \bigr)_{q \in \kx}$ be a flat family of quantized coordinate rings for an affine variety $V$, where $\A$ is again a torsionfree algebra over a Laurent polynomial ring $k[t^{\pm1}]$ and $A_1 \cong \calO(V)$. 

Further, assume that $\A$ is defined over a subfield $k_0 \subset k$, that is, $\A = \A_0 \otimes_{k_0} k$ for a torsionfree $k_0[t^{\pm1}]$-algebra $\A_0$. 
We claim that if $p,q \in \kx$ are transcendental over $k_0$, then $\spec A_p \approx \spec A_q$ and $\prim A_p \approx \prim A_q$.

To see this, observe that there is a $k_0$-automorphism $\phi \in \Aut k$ such that $\phi(p) = q$. Then $\phi$ extends to a $k_0[t^{\pm1}]$-algebra automorphism $\phihat = \id_{\A_0} \otimes\, \phi$ of $\A$, and $\phihat$ induces an isomorphism $A_p \longrightarrow A_q$ of $k_0$-algebras. This isomorphism, in turn, induces the claimed homeomorphisms.

Let us consider these homeomorphisms in a particular instance, namely when $\A := k[t^{\pm1}] \langle x,y \mid xy = t yx \rangle$. This algebra determines a flat family $\bigl( \Oq(k^2) \bigr)_{q\in\kx}$ and is defined over the prime field $k_0$ of $k$. Consequently, the quantum planes $\calO_p(k^2)$ and $\Oq(k^2)$ are isomorphic as rings for any $p,q \in \kx$ which are transcendental over $k_0$. On the other hand, these rings are not isomorphic as $k$-algebras unless $p = q^{\pm1}$, as follows from \cite[Proposition 3.2, Corollaire 3.11(c)]{JAFD}.

On the side, we point out that the homeomorphisms arising from $\phihat$ in this instance are not entirely geometric, in that their restrictions to classical portions of the given spectra are not always morphisms of classical varieties.  The homeomorphism $\prim A_p \longrightarrow \prim A_q$ induced by $\phihat$ sends
$$\langle x-\alpha,\, y\rangle A_p \ \ \longmapsto \ \ \langle x- \phi(\alpha),\, y\rangle A_q$$
for $\alpha \in \kx$. In other words, $\kx$ is being mapped to itself via $\phi$, and the latter map is not a morphism of $k$-varieties unless $p=q$ and $\phi = \id_k$.
\end{observation}

\section{Piecewise classical structure of spectra}  \label{piecewise}
The prime and primitive spectra of many quantum algebras, particularly generic quantized coordinate rings, are known to have a ``piecewise classical" structure, in that they possess finite stratifications in which each stratum is homeomorphic to a classical affine scheme or variety. We outline the relevant partitions, raise the question of fully describing these spectra in terms of classical data, and set up a conjectural framework for achieving an integrated picture.

\begin{theorem} \label{stratif}  \textup{ \cite[Theorems II.2.13, II.8.4, Proposition II.8.3]{KBKG02}}
Let $A$ be a noetherian $k$-algebra and 
$H$ a torus $(\kx)^r$ acting rationally on $A$. 

\textup{(a)} There are partitions
$$\spec A = \bigsqcup_{J \in \Hspec A} \spec_J A\quad\text{and}\quad \prim A = \bigsqcup_{J \in \Hspec A} \prim_J A,$$
where $\Hspec A$ is the set of $H$-stable prime ideals of $A$ and
$$\spec_J A := \{ P \in \spec A \mid \bigcap_{h\in H} h(P) = J \},\qquad \prim_J A := \prim A \cap \spec_J A.$$

\textup{(b)} For each $J \in \Hspec A$, there are homeomorphisms $\spec_J A \approx \spec Z_J$ and $\prim_J A \approx \max Z_J$, where $Z_J$ is the center of a localization of $A/J$ and is isomorphic to a Laurent polynomial ring over $k$.

\textup{(c)} Now assume additionally that $\Hspec A$ is finite, and that $A$ satisfies the noncommutative Nullstellensatz over $k$. Then
$$\prim_J A = \{\,\text{\rm maximal elements of}\; \spec_J A\,\}$$
for all $J \in \Hspec A$.
\end{theorem}

The above partition of $\spec A$ is a \emph{stratification} in the sense that the closure of any stratum $\spec_J A$ is a union of strata, and likewise for the partition of $\prim A$. For details on quantized coordinate rings that satisfy the hypotheses of Theorem \ref{stratif}, see \cite[Corollary II.8.5 and proof]{KBKG02}.

The topologies on these spectra are determined by the topologies on the strata together with relations between the strata, such as conditions on relative closures. In order to describe these topologies, we write $\Cl(T)$ for the collection of closed subsets of a topological space $T$, and we use overbars to denote closures in $\spec A$ and $\prim A$.

\begin{lemma}  \label{descripclosed}  \textup{\cite[Lemma 2.3, Observation 3.4]{KBKG15}}
Let $A$ and $H$ be as in Theorem \textup{\ref{stratif}}.
For $J \subset J'$ in $\Hspec A$, define
\begin{align*}
\phi^s_{JJ'} &: \Cl(\spec_J A) \rightarrow \Cl(\spec_{J'} A), & \phi^s_{JJ'}(Y) &= \overline{Y}\cap \spec_{J'} A  \\
\phi^p_{JJ'} &: \Cl(\prim_J A) \rightarrow \Cl(\prim_{J'} A), & \phi^p_{JJ'}(Y) &= \overline{Y}\cap \prim_{J'} A\,.
\end{align*}

Then the closed subsets of $\spec A$ are the subsets $X$ such that
\begin{itemize}
\item $X \cap \spec_J A \in \Cl(\spec_J A)$ for all $J$;
\item $\phi^s_{JJ'}(X \cap \spec_J A) \subseteq X \cap \spec_{J'} A$ for all $J \subseteq J'$.
\end{itemize}
Similarly, the closed subsets of $\prim A$ are the subsets $X$ such that
\begin{itemize}
\item $X \cap \prim_J A \in \Cl(\prim_J A)$ for all $J$;
\item $\phi^p_{JJ'}(X \cap \prim_J A) \subseteq X \cap \prim_{J'} A$ for all $J \subseteq J'$.
\end{itemize}
\end{lemma}

Since the strata of $\spec A$, respectively $\prim A$, are classical schemes, respectively varieties, one would like to also describe the maps $\phi^\bullet_{JJ'}$ in terms of classical geometric data. The most obvious possibilities fail, however. Immediate examples show that these  maps cannot always be given by inverse images under morphisms, nor by closures of images under morphisms (cf.~\cite[Example 3.5]{KBKG15} and Observation \ref{auxOqGL2} below). Potentially, the $\phi^\bullet_{JJ'}$ can be described in terms of morphisms connecting the relevant spectra to auxiliary schemes or varieties as follows, where we identify $\prim_J A$ and $\spec_J A$ with the spaces $\max Z_J$ and $\spec Z_J$ given in Theorem \ref{stratif}(b).

\begin{conjecture} \label{closuremaps} \textup{\cite[Conjecture 3.11]{KBKG15}} 
Let $A$ and $H$ be as in Theorem \textup{\ref{stratif}(c)}.

For $J \subset J'$ in $\Hspec A$, there exist an affine variety $V_{JJ'}$, a corresponding scheme $\Vtil_{JJ'}$, and morphisms
$$\xymatrixrowsep{3pc} \xymatrixcolsep{6pc}
\xymatrix{
\prim_{J'} A \ar[dr]^{f_{JJ'}} &&\spec_{J'} A \ar[dr]^{\ftil_{JJ'}}  \\
\prim_J A \ar[r]^{g_{JJ'}} &V_{JJ'} &\spec_J A \ar[r]^{\gtil_{JJ'}} &\Vtil_{JJ'}
}$$
such that
\begin{enumerate}
\item[] $\phi^p_{JJ'}(Y) = f_{JJ'}^{-1}(\,\overline{g_{JJ'}(Y)}\,)$ for $Y \in \Cl(\prim_J A)$;
\item[] $\phi^s_{JJ'}(Z) = \vphantom{3^{3^{3^{3^J}}}} 
\ftil_{JJ'}^{-1}(\,\overline{\gtil_{JJ'}(Z)}\,)$ for $Z \in \Cl(\spec_J A)$.
\end{enumerate}
\end{conjecture} 

Conjecture \ref{closuremaps} was confirmed for $\Oq(GL_2(k))$, $\Oq(SL_3(k))$ and $\Oq(M_2(k))$ in \cite[Example 4.3, Theorems 5.4, 7.5]{KBKG15}; Poisson analogs for $\calO(GL_2(k))$ and $\calO(SL_3(k))$ were established in \cite[Corollaries 4.5, 4.6]{Fry}.

\begin{remark}  \label{closuremaps.implies.semiclass}
Suppose that $k$, $A$, $V$ are as in the hypotheses of Conjecture \ref{semiclassconj}, and that the hypotheses of Theorem \ref{stratif}(c) hold for a suitable torus $H$. In view of Theorem \ref{stratif} and Lemma \ref{descripclosed}, we see that if Conjecture \ref{closuremaps} is validated, then the spaces $\spec A$ and $\prim A$ are determined by the strata $\spec_J A$ and $\prim_J A$ together with the auxiliary data from the latter conjecture.

This leads to further bridges among the conjectures.  Indeed, if one assumes  Conjecture \ref{closuremaps} as well as a matching Poisson analogue for  $\calO(V)$, one obtains homeomorphisms $\Pspec \calO(V) \approx \spec A$ and $\Pprim \calO(V) \approx \prim A$. This approach was used to verify Conjecture \ref{semiclassconj} for $\Oq(SL_3(k))$ and $\OqGLtwok$ in \cite[Theorem 5.21, Corollaries 5.22, 5.23]{Fry}.
\end{remark}

\begin{example}  \label{auxOqGL2}
To illustrate the auxiliary varieties and morphisms that may appear in the context of Conjecture \ref{closuremaps}, consider the algebra $A = \OqGLtwok$, $q\nrootone$, and the auxiliary data that have been used to confirm the conjecture in this case \cite[Example 4.3]{KBKG15}.

There is a natural action of the torus $H = (\kx)^4$ on $A$, under which $A$ has exactly $4$ $H$-primes; the poset $\Hspec A$ is drawn on the left in Figure \ref{HspecOqGL2}. The strata $\prim_J A$ are homeomorphic to the varieties shown on the right. We identify the $\prim_J A$ with these varieties.
\begin{figure}[ht]
$$\xymatrixrowsep{1pc} \xymatrixcolsep{1pc}
\xymatrix{
&\langle b,c\rangle &&&&&(\kx)^2  \\
\langle b\rangle \edge[ur] &&\langle c\rangle \edge[ul] &&&\kx &&\kx  \\
&\langle 0\rangle \edge[ul] \edge[ur] &&&&&(\kx)^2
}$$
\caption{$\Hspec \OqGLtwok$ and strata $\prim_J \OqGLtwok$}
\label{HspecOqGL2}
\end{figure}

The auxiliary data used for $\prim A$ \cite[Figure 1]{KBKG15} are displayed in Figure \ref{auxprimOqGL2}, where the varieties $V_{JJ'}$ are marked with square brackets. (The data for $\phi^p_{\langle b\rangle,\langle b,c\rangle}$ and $\phi^p_{\langle c\rangle,\langle b,c\rangle}$ are the same, so we display only one set.)
\begin{figure}[ht]
$$\xymatrixrowsep{4pc} \xymatrixcolsep{6pc}
\xymatrix{
&(\kx)^2 \ar[ddl]_{\operatorname{mult}} \ar[dr]^{\operatorname{mult}}  \\
&\kx \ar[r]^{\operatorname{id}} \ar[dr]^{(0,-)} &\bigl[\,\kx\,\bigr] \\
\bigl[\,\kx\,\bigr] &(\kx)^2 \ar[l]_{\operatorname{pr}_2} \ar[r]^{\operatorname{incl}} &\bigl[\,k\times\kx\,\bigr] 
}$$
\caption{Auxiliary data for $\prim \OqGLtwok$}
\label{auxprimOqGL2}
\end{figure}

Obviously the map $\phi^p_{\langle b\rangle,\langle b,c\rangle}$ can be expressed in the form $Y \mapsto h_b^{-1}(Y)$, where $h_b : (\kx)^2 \rightarrow \kx$ is the multiplication map. Describing the maps $\phi^p_{\langle0\rangle, \langle b\rangle}$ and $\phi^p_{\langle0\rangle,\langle b,c\rangle}$ in similar terms is a different story, however. For instance, $\phi^p_{\langle0\rangle, \langle b\rangle}$ sends all finite subsets of $(\kx)^2$ to the empty set. Hence, it cannot be given by inverse images under any map $\kx \rightarrow (\kx)^2$, nor can it be given by closures of images under any map $(\kx)^2 \rightarrow \kx$.
\end{example}

\end{document}